\documentclass[12pt]{amsart}

\usepackage{amssymb}
\usepackage{amscd}

\newtheorem{theorem}{Theorem}[section]
\newtheorem{proposition}[theorem]{Proposition}
\newtheorem{corollary}[theorem]{Corollary}

\theoremstyle{remark}
\newtheorem{example}[theorem]{Example}

\numberwithin{equation}{section}

\font\tenscr=rsfs10 at 12pt

\newcommand{\Ga}{\mathbf{G}_{\mathrm{a}}}
\newcommand{\Gm}{\mathbf{G}_{\mathrm{m}}}
\newcommand{\bG}{\mathbf{G}}
\newcommand{\kernel}{\mathop{\mathrm{Ker}}\nolimits}
\newcommand{\Spec}{\mathop{\mathrm{Spec}}\nolimits}

\begin{document}

\title[On a discrepancy among Picard-Vessiot theories]{On a discrepancy among Picard-Vessiot
theories in positive characteristics}
\author[Katsutoshi Amano]{Katsutoshi Amano}
\address{8-29, Ida-Sammai-cho, Nakahara, Kawasaki, Kanagawa, 211-0037, Japan}
\email{ma-pfybdb-612019@agate.dti.ne.jp}
\address{{\rm http://www.green.dti.ne.jp/amano/index-eng.html}}
\date{Sep.\ 2007; revised version 2.}

\begin{abstract}
There is a serious discrepancy among literature on the Picard-Vessiot theory in positive characteristics (for iterative differential fields).
We consider descriptions of Galois correspondence in four approaches to this subject: Okugawa's result \cite{O}, Takeuchi's Hopf algebraic approach \cite{T} (and \cite{AM}), the result of Matzat and van der Put \cite{MP}, and the model theoretic approach by Pillay \cite{P}.
In the three approaches except Takeuchi's one, Galois correspondence is described between closed subgroups of algebraic matrix groups and its fixed fields.
But such a description has a problem that the Galois correspondence may not be bijective there.
We explain this problem in the first section by giving an explicit example.

We should use affine group schemes instead of algebraic matrix groups to obtain a suitable Galois correspondence, as in Takeuchi's approach.
But intermediate fields are not necessarily fixed fields there.
In the second section, we give some sufficient conditions for intermediate artinian simple module algebras to be fixed algebras in the context of the unified Picard-Vessiot theory developed by the author and Masuoka \cite{AM}.
\end{abstract}

\maketitle

\baselineskip 18pt

\section*{Introduction}

We compare descriptions of Galois correspondence in four approaches to the Picard-Vessiot theory in positive characteristics (for iterative differential fields): Okugawa's result \cite{O}, Takeuchi's approach \cite{T} (and \cite{AM}), the result of Matzat and van der Put \cite{MP}, and the model theoretic approach by Pillay \cite{P}.
All descriptions except Takeuchi's approach are written between closed subgroups of Galois group and its fixed fields.
But, in positive characteristics, such a description has a problem that the Galois correspondence may not be bijective.
(Actually, Okugawa avoided this problem by restricting his definition of Picard-Vessiot extension.)
In Section \ref{sec-example}, we discuss this point by considering an explicit example.
Especially we emphasize certain advantages of Takeuchi's approach: a complete description of bijective Galois correspondence is obtained only by his approach.

There is another interest came from a more innocent motivation.
In \cite{AM}, the author and Masuoka developed a generalized and unified Picard-Vessiot theory including Galois theories for differential modules, difference modules, or mixed cases of them in arbitrary characteristics.
Our description of Galois correspondence is based on Takeuchi's approach: the correspondence among intermediate artinian simple (AS) $D$-module algebras of a Picard-Vessiot extension $L/K$, $D$-stable coideals of $L\otimes_{K}L$, and Hopf ideals.
But intermediate AS $D$-module algebras may not be ``fixed algebras" in the usual sense.
The second interest of this article is to investigate this problem in the generalized context of \cite{AM}.
In Section \ref{sec-criteria}, we give some criteria to see which intermediate AS $D$-module algebras are ``fixed algebras" in the usual sense. 
Also, by this result, we can have a clear insight into our first problem on non-bijectivity of Galois correspondence in terms of algebraic matrix groups.

To describe main results in Section \ref{sec-criteria}, we use terminologies of \cite{AM} (see also \cite[Part 3]{Thesis}).
Here we outline it briefly. We assume that the reader is familiar with concepts about Hopf algebras and affine group schemes.
Let $k$ be a field and $D$ a cocommutative pointed Hopf algebra over $k$ whose irreducible component $D^{1}$ containing $1$ is a Birkhoff-Witt bialgebra.
A $D$-module algebra is called {\em artinian simple} (AS) iff it is artinian as a ring and it has no nontrivial $D$-stable ideal.
An extension of AS $D$-module algebras $L/K$ is called a {\em Picard-Vessiot extension} iff $K^{D}=L^{D}$ and there exists a (necessarily unique) $D$-module subalgebra $K\subset A\subset L$ such that the total quotient ring of $A$ equals $L$ and $H=(A\otimes_{K}A)^{D}$ generates the left $A$-module $A\otimes_{K}A$.
Then $H$ has a structure of commutative Hopf algebra over $K^{D}$.
We say that $(L/K,A,H)$ is a Picard-Vessiot extension to indicate these.
There is a bijective correspondence among intermediate AS $D$-module algebras $F$ of $L/K$, $D$-stable coideals $J$ of the $L$-coring $L\otimes_{K}L$, and Hopf ideals $I$ of $H$ as follows:
\[ \begin{array}{l}
F \mapsto J=\kernel(L\otimes_{K}L\twoheadrightarrow L\otimes_{F}L)
\mapsto I=H\cap J, \vspace{5pt} \\
I \mapsto J=I\cdot(L\otimes_{K}L)\mapsto F=\{a\in L\;|\; 1\otimes a-a\otimes 1\in J\}.
\end{array} \]
If $F$ is an intermediate AS $D$-module algebra corresponding to a Hopf ideal $I\subset H$, then $(L/F,AF,H/I)$ is also a Picard-Vessiot extension.
The affine group scheme $\bG(L/K):=\Spec H$ represented by $H$ is called the {\em Picard-Vessiot group scheme} for $L/K$.
This is isomorphic to the automorphism group functor $\mathbf{Aut}_{D}(A/K)$ which associates each commutative $K^{D}$-algebra $T$ the $D$-linear $K\otimes_{K^{D}}T$-algebra automorphism group:
\[ \bG(L/K)\simeq\mathbf{Aut}_{D}(A/K) : T\mapsto\mathrm{Aut}_{D}(A\otimes_{K^{D}}T/K\otimes_{K^{D}}T). \]
Especially $\bG(L/K)(K^{D})$ is isomorphic to $\mathrm{Aut}_{D}(A/K)=\mathrm{Aut}_{D}(L/K)$.
$L/K$ is called {\em finitely generated} iff there exist finite $x_{1},\dotsc,x_{n}\in L$ such that $L=K\langle x_{1},\dotsc,x_{n}\rangle$, the smallest AS $D$-module subalgebra of $L$ including both $K$ and $x_{1},\dotsc,x_{n}$.
This is the case iff $A$ is a finitely generated $K$-algebra iff $H$ is a finitely generated $K^{D}$-algebra.

Let $(L/K,A,H)$ be a Picard-Vessiot extension of AS $D$-module algebras such that $L^{D}=K^{D}=k$.
As our main results, we show the following theorems. \\

\noindent
{\bf Theorem \ref{thm-FF}}. {\it Assume that $L/K$ is finitely generated.
Let $F$ be an intermediate AS $D$-module algebra of $L/K$ corresponding to a Hopf ideal $I\subset H$.
If $H/I$ satisfies
\[ \bigcap_{g\in\mathrm{Alg}_{k}(H/I,k)}\kernel g=0, \]
then $L^{\bG(L/F)(k)}=F$.
} \vspace{5pt}

We observe the last sufficient condition of the above theorem is equivalent to saying that $H/I$ coincides with the coordinate ring of the affine algebraic group $\bG(L/F)(k)$; see \cite[(4.5)]{W} and the remark after (\ref{AAG}) in this article.
When $k$ is algebraically closed, it is equivalent to saying that $H/I$ is reduced.
In Proposition \ref{prop-separable}, we will give some equivalent conditions to this.

As a generalization of the Galois correspondence in Okugawa's sense \cite[Ch.\ II, \S 7]{O}, we will show: \\

\noindent
{\bf Theorem \ref{thm-Okugawa}}. {\it Assume that $k$ is algebraically closed, $H$ is reduced, and $L/K$ is finitely generated.
Let $\mathcal{F}$ be the set of intermediate AS $D$-module algebras $F$ of $L/K$ such that $AF\otimes_{F}AF$ is reduced.
(When $L$ is a field, $\mathcal{F}$ becomes the set of intermediate $D$-module fields of $L/K$ over which $L$ is (classically) separable.)
Let $\mathcal{G}$ be the set of all closed subgroups of $\bG(L/K)(k)$.
Then $\mathcal{F}$ and $\mathcal{G}$ corresponds bijectively by
\[ \begin{array}{l}
\mathcal{F}\rightarrow\mathcal{G},\quad F\mapsto\bG(L/F)(k),
\vspace{5pt} \\
\mathcal{G}\rightarrow\mathcal{F},\quad \mathcal{H}\mapsto L^{\mathcal{H}}.
\end{array} \]
} \vspace{5pt}

\section{An example of inseparable Picard-Vessiot extension of iterative differential fields}
\label{sec-example}

We consider fields with higher derivations $\{\partial^{(n)}\}_{n\geq 0}$ of infinite length.
In \cite{MP}, such fields are called {\em iterative differential (ID) fields}; we follow this terminology for a while.
There are four approaches to the Picard-Vessiot theory for them:
\begin{enumerate}
\renewcommand{\labelenumi}{(\arabic{enumi})}
\item Okugawa's result \cite{O}; \label{Okugawa}
\item Takeuchi's Hopf algebraic approach \cite{T}
(see also \cite{AM}, \cite[Part 3]{Thesis}); \label{Takeuchi}
\item The result of Matzat and van der Put \cite{MP}; \label{Matzat}
\item The model theoretic approach by Pillay \cite{P}. \label{Pillay}
\end{enumerate}
We make a comparison among them, especially on Galois correspondence.

It should be noted that Okugawa's definition of Picard-Vessiot extension in \cite{O} is properly stronger than others (see Kolchin's review on \cite{O} at MathSciNet) and a bijective Galois correspondence is obtained.
But existence of Picard-Vessiot extensions for ID-modules cannot be shown in terms of (\ref{Okugawa}).
The definitions of Picard-Vessiot extensions of others are essentially equivalent (see \cite[Theorem 3.3]{T} or \cite[Theorem 4.6]{AM}) and existence of Picard-Vessiot extensions is shown in each approach.
But descriptions of Galois correspondence are different.
In (\ref{Matzat}) and (\ref{Pillay}), the correspondence is described between closed subgroups of algebraic matrix groups and its fixed fields.
In (\ref{Takeuchi}), algebraic matrix groups are replaced with affine group schemes (commutative Hopf algebras) and the correspondence is described among closed subgroup schemes (Hopf ideals), $C$-ferential coideals, and intermediate $C$-ferential fields.
Actually, there is a little mistake in the proof of Galois correspondence in (\ref{Matzat}) and (\ref{Pillay}), and their Galois correspondence is not bijective.
Only (\ref{Takeuchi}) succeeds to prove both of existence of Picard-Vessiot extensions and a bijective Galois correspondence, however intermediate fields are no longer ``fixed fields".
We will see these facts by the following example.

As a matter of fact, the author found this example by following \cite[p.\ 14]{MP}.
First we construct a Picard-Vessiot extension whose Galois group is $\Ga$.
Let $k=\bar{\mathbb{F}}_{p}$ be the algebraic closure of the prime field in characteristic $p>0$.
We define an ID-ring structure on the polynomial ring $k[z]$ by
\[ \partial^{(n)}z^{m}=\left(\begin{array}{c}m \\ n \end{array}\right)z^{m-n} \]
with $k$ the field of constants.
This structure is uniquely extended to the quotient field $K=k(z)$ and $K$ becomes an ID-field.
Let $V=Kv_{1}+Kv_{2}$ be a two-dimensional $K$-vector space with a basis $v_{1}, v_{2}$.
Define an ID-module structure on $V$ by
\[ \partial^{(n)}v_{2}=0,\quad \partial^{(n)}v_{1}=\left\{\begin{array}{ll} v_{2} & (n=p^{p^{l}}\ \mbox{for some}\ l\in\mathbb{Z}_{\geq 0}) \vspace{5pt} \\ 0 & \mbox{(otherwise),} \end{array}\right. \]
for $n=1,2,3,\dotsc$.
Then a solution of $V$ is given by the following formal power series:
\[ f:=\sum_{l=0}^{\infty}z^{p^{p^{l}}}\in k[[z]]. \]
By \cite[Lemma 5.2]{MP}, we have that $f$ is transcendental over $K$.
Put $L=K(f)$. Then $L/K$ is a Picard-Vessiot extension for $V$ and the Galois group (in their sense) is $G(L/K)=\Ga(k)$.
The action of $\sigma\in k=\Ga(k)$ on $L$ is given by $f\mapsto f+\sigma$.

Then let $F=K(f^{p})$ ($\subsetneq L$).
This is an intermediate ID-field of $L/K$.
But there is no closed subgroup $H\subset\Ga(k)$ such that $L^{H}=F$.
This implies that the correspondence described in \cite[Theorem 3.5]{MP} is not surjective.
(The proof corresponding to \cite[Proposition 3.12 (3)]{P2} fails since $R\otimes_{K}R$ is not necessarily reduced.)

In Pillay's proof of Galois correspondence \cite[Proposition 3.2, p.\ 334]{P}, it seems that \cite[Fact 3.2, p.\ 335]{P}, which asserts that every Picard-Vessiot extension of ID-fields is a (classically) separable extension, is needed.
But $L/F$ is a counter example since it is a Picard-Vessiot extension for $F\otimes_{K}V$ and is purely inseparable. 

We observe that $L/F$ is not a Picard-Vessiot extension in Okugawa's sense.
His definition of Picard-Vessiot extension is properly stronger than ours.
He called an ID-field extension $M/N$ Picard-Vessiot if it satisfies the following condition (S) in addition to our definition: \\

\noindent
{\bf Condition (S)}.
{\it For every $x\in M-N$, there exists an ID-field extension $M'\supset M$ and a differential $N$-algebra map $\sigma : M\rightarrow M'$ such that $\sigma(x)\neq x$.} \\

\noindent
If this condition satisfied, we say that $M$ satisfies (S) over $N$.
We can show directly that $L$ does not satisfy (S) over $F$. In fact, for any ID-field extension $L'\supset L$ and differential $F$-algebra map $\sigma : L\rightarrow L'$, we have $\sigma(f)^{p}=\sigma(f^{p})=f^{p}$ and hence $\sigma(f)=f$.
According to Kolchin's review, the condition (S) is equivalent to saying that the extension is (classically) separable, provided other conditions in the definition of Picard-Vessiot extension.
(In Section \ref{sec-criteria}, we will verify this fact implicitly.)
Though $L/K$ is a Picard-Vessiot extension also in Okugawa's sense, $F$ does not appear in his description of Galois correspondence \cite[Ch.\ II, \S 7]{O}. 
Let $\mathcal{F}$ be the set of intermediate ID-fields of $L/K$ over which $L$ satisfies the condition (S) and $\mathcal{G}$ the set of all closed subgroups of $G(L/K)$.
Then Okugawa's description of Galois correspondence is as follows:
\[ \begin{array}{l}
\mathcal{F}\rightarrow\mathcal{G},\quad  M\mapsto G(L/M),
\vspace{5pt} \\
\mathcal{G}\rightarrow\mathcal{F},\quad  H\mapsto L^{H},
\end{array} \]
which is bijective.
We will generalize this correspondence in Theorem \ref{thm-Okugawa}.

In the following, we describe the Picard-Vessiot extension $L/K$ in terms of \cite{AM} and explain that $F$ is caught only by an affine group scheme in Galois correspondence.
Let $D=B(k)=\sum_{n=0}^{\infty}k\partial^{(n)}$ be the Birkhoff-Witt bialgebra spanned by one divided power sequence $\{\partial^{(n)}\}$ of infinite length (see e.g.\ \cite[Ch.\ 2, \S 5]{Abe}) over $k$.
Then $D$ is a cocommutative irreducible Hopf algebra.
The coalgebra structure $(\Delta,\varepsilon)$ and the algebra structure is given by
\[ \Delta\partial^{(n)}=\sum_{i+j=n}\partial^{(i)}\otimes\partial^{(j)},\ \ \varepsilon(\partial^{(n)})=\left\{\begin{array}{ll} 1 & (n=0) \vspace{5pt} \\ 0 & (n>0), \end{array}\right. \quad\partial^{(n)}\partial^{(m)}=\left(\begin{array}{c} n+m \\ n \end{array}\right)\partial^{(n+m)}. \]
It is known that an irreducible bialgebra necessarily have the antipode \cite[Theorem 2.4.24]{Abe}; the antipode $S$ of $D$ is given by $S(\partial^{(n)})=(-1)^{n}\partial^{(n)}$.
In this case the notion of $D$-module algebras is the same as the notion of ID-rings whose constants include $k$.
The category of ID-modules over $K$ is the same as the category of $K\# D$-modules ${}_{K\#D}\mathcal{M}$.
We have that $L/K$ is a minimal splitting field for $V$ and hence it is a Picard-Vessiot extension \cite[Theorem 4.6]{AM}.
The Picard-Vessiot group scheme for $L/K$ is $\bG(L/K)=\Ga$ and the closed subgroup scheme corresponding to $F$ is $\boldsymbol{\alpha}_{p}\subset\Ga$.
Since $\boldsymbol{\alpha}_{p}(k)$ is trivial, it does not appear as a nontrivial closed subgroup of $\Ga(k)$.
Thus $F$ cannot be caught as a fixed field of an algebraic matrix group.

The principal algebra (or the Picard-Vessiot ring) for $L/K$ is $A=K[f]$ and the Hopf algebra representing $\bG(L/K)$ is $H=(A\otimes_{K}A)^{D}=k[l]$ with one primitive $l=1\otimes f-f\otimes 1$.
Every Hopf ideal of $H$ is written as $\langle\varphi(l)\rangle$ where $\varphi(x)=\sum_{i}a_{i}x^{p^{i}}\in k[x]$; see \cite[Ch.\ 8, (8.4) and Ex.\ 7]{W}.
The Galois correspondence for $L/K$ in terms of \cite{AM} is as follows:
\[ \begin{array}{ccccc}
\mbox{\{intermediate ID-field\}} & &
\mbox{\{$D$-stable coideal of $L\otimes_{K}L$\}}
 & & \mbox{\{Hopf ideal of $H$\}} \vspace{5pt} \\
K(\varphi(f)) & \leftrightarrow & L\cdot\varphi(l)\cdot L
 & \leftrightarrow & \langle \varphi(l)\rangle. \end{array} \]
If $a_{0}\neq 0$, $L/K(\varphi(f))$ is a finite separable extension and $\bG(L/K(\varphi(f)))=\Spec H/\langle\varphi(l)\rangle$ is finite etale.
Otherwise $L/K(\varphi(f))$ is inseparable, $H/\langle\varphi(l)\rangle$ is not reduced, and the fixed field by $\bG(L/K(\varphi(f)))(k)$ is not $K(\varphi(f))$ (it is a purely inseparable extension of $K(\varphi(f))$).

We can also describe the above argument in terms of \cite{T} by taking $C=D$ which is considered merely as an irreducible coalgebra.
We have $T(C^{+})=k\langle \partial^{(1)},\partial^{(2)},\dotsc\rangle$, the (non-commutative) free algebra generated by $\partial^{(1)}, \partial^{(2)},\dotsc$.
This also becomes a Birkhoff-Witt bialgebra since the $\hbox{\tenscr V}$-map preserves multiplications (see \cite{Heyneman-Sweedler}).
In this case the category of ID-modules over $K$ is a full abelian tensor subcategory of the category ${}_{K\#T(C^{+})} \mathcal{M}$ of $C$-ferential $K$-modules.

\section{When an intermediate AS module algebra is a fixed algebra?}
\label{sec-criteria}

Let $k$ be an arbitrary field and $\bG$ an algebraic affine $k$-group scheme.
Consider the following condition on the coordinate Hopf algebra $H=k[\bG]$ representing $\bG$:
\begin{equation}
\label{AAG}
\bigcap_{g\in\mathrm{Alg}_{k}(H,k)}\kernel g=0.
\end{equation}
By \cite[(4.5)]{W}, this is equivalent to saying that $H$ coincides with the coordinate ring of the affine algebraic group $\bG(k)$; namely, (\ref{AAG}) $\Leftrightarrow$ $H=k[\bG(k)]$.
(Though it is assumed that $k$ is an infinite field in \cite[Ch.\ 4]{W}, Theorem in \cite[(4.5)]{W} holds for arbitrary base fields.
Essentially, the assumption is needed only for \cite[(4.5), Corollary, p.\ 32]{W} and ``only if" part of the proof of \cite[(4.6)]{W}.)
When $k$ is algebraically closed, we have that (\ref{AAG}) $\Leftrightarrow$ $H$ is reduced $\Leftrightarrow$ $\bG$ is smooth.
Since the common kernel in (\ref{AAG}) includes the nilradical of $H$, (\ref{AAG}) implies that $H$ is reduced.
But the converse does not hold in general; see the following example.

\begin{example}
(1) $k[\Ga]=k[l]$ (resp., $k[\Gm]=k[x,x^{-1}]$) satisfies (\ref{AAG}) iff $k$ is an infinite field.
When $k$ is a finite field of order $q$, the common kernel is $\langle l^{q}-l\rangle$ (resp., $\langle x^{q-1}-1\rangle$).

(2) Let $\bG$ be finite etale. Let $k_{s}$ be the separable closure of $k$.
Then $\mathrm{Gal}(k_{s}/k)$ acts on $\bG(k_{s})$ as group automorphisms.
The action is trivial iff $k[\bG]$ satisfies (\ref{AAG}).
See \cite[Example 3.7.6]{Thesis} for an example of a Picard-Vessiot group scheme which is finite etale but does not satisfy (\ref{AAG}).
\end{example}

Let $D$ be a cocommutative, pointed $k$-Hopf algebra whose irreducible component $D^{1}$ containing $1$ is a Birkhoff-Witt bialgebra, as in \cite{AM}.
Let $(L/K,A,H)$ be a finitely generated Picard-Vessiot extension of AS $D$-module algebras such that $L^{D}=K^{D}=k$.
Our first criterion is obtained as follows:
\begin{theorem}
\label{thm-FF}
$L^{\bG(L/K)(k)}=K$ if $H$ satisfies (\ref{AAG}).
\end{theorem}
\begin{proof}
This follows in the same way as \cite[Proposition 3]{Levelt}.
Let $\mu : A\otimes_{k}H\rightarrow A\otimes_{K}A$, $a\otimes h\mapsto ah$ be the isomorphism in \cite[Proposition 3.4]{AM}.
Recall that the $H$-comodule structure of $A$ is given by $A\rightarrow A\otimes_{k}H$, $a\mapsto \mu^{-1}(1\otimes a)$.

Let $x$ be an arbitrary element in $L^{\bG(L/K)(k)}$.
Write $x=b/c$ with $b,c\in A$.
Put $w=\mu^{-1}(b\otimes c-c\otimes b)\in A\otimes_{k}H$.
For $g\in\bG(L/K)(k)=\mathrm{Alg}_{k}(H,k)$, the action of $g$ on $x\in L$ is given by $\rho_{g}(x):=(\mathrm{id}_{A}\otimes g)(\mu^{-1}(1\otimes b))/(\mathrm{id}_{A}\otimes g)(\mu^{-1}(1\otimes c))$.
Since $\rho_{g}(x)=x$ for all $g\in\mathrm{Alg}_{k}(H,k)$ and since $x=(\mathrm{id}_{A}\otimes g)(\mu^{-1}(b\otimes 1))/(\mathrm{id}_{A}\otimes g)(\mu^{-1}(c\otimes 1))$, we have $(\mathrm{id}_{A}\otimes g)(w)=0$ for all $g\in\mathrm{Alg}_{k}(H,k)$.
We claim $w=0$.
To see this, write $w=\sum a_{i}\otimes h_{i}$, where $a_{i}\in A$ are $k$-linearly independent and $h_{i}\in H$.
Then $0=(\mathrm{id}_{A}\otimes g)(w)=\sum a_{i}g(h_{i})$, and hence $g(h_{i})=0$ for all $i$ and for all $g\in\mathrm{Alg}_{k}(H,k)$.
By the condition (\ref{AAG}), we have $h_{i}=0$ for all $i$.
Therefore $w=0$ and $b\otimes c=c\otimes b$, which implies $x\in K$.
\end{proof}

\begin{corollary}
\label{cor}
Let $F$ be an intermediate AS $D$-module algebra of $L/K$, corresponding to a Hopf ideal $I\subset H$.
Then $L^{\bG(L/F)(k)}=F$ if $H/I$ satisfies (\ref{AAG}).
\end{corollary}

Let $G=G(D)$ be the group of grouplikes in $D$.
Fix a maximal ideal $P$ of $K$ and put $G_{P}=\{g\in G\;|\;gP=P\}$.
Recall that $K$ is isomorphic to a product of some finite copies of a field $K_{1}=\Psi_{G_{P}}(K)$ \cite[Proposition 2.4]{AM}.
Let $L_{1}=\Psi_{G_{P}}(L)$ and $A_{1}=\Psi_{G_{P}}(A)$.
On the condition (\ref{AAG}), we have:
\begin{proposition}
\label{prop-separable}
When $k$ is algebraically closed, the following are equivalent:
\begin{enumerate}
\renewcommand{\labelenumi}{(\alph{enumi})}
\item $H$ is reduced ($\Leftrightarrow$ $\bG(L/K)$ is smooth
$\Leftrightarrow$ $H$ satisfies (\ref{AAG}));
\item $A_{1}$ is a (classically) separable $K_{1}$-algebra;
\item $L_{1}$ is a (classically) separable $K_{1}$-algebra;
\item $A\otimes_{K}A$ is reduced.
\end{enumerate}
\end{proposition}
\begin{proof}
Since $L$ is reduced, $A$ is also reduced. Thus (a) $\Leftrightarrow$ (d) follows from the isomorphism $\mu : A\otimes_{k}H\rightarrow A\otimes_{K}A$ (see also \cite[Ch.\ 6, Ex.\ 2]{W}).
To see (a) $\Leftrightarrow$ (b), we may assume that $K$ is a field since $(L_{1}/K_{1},A_{1},H)$ is a Picard-Vessiot extension of AS $D(G_{P})$-module algebras \cite[Lemma 3.7]{AM}.
Let $M\subset A$ be a maximal ideal of $A$.
Since $A$ is a finitely generated $K$-algebra, we have $A/M$ is isomorphic to a subfield of $\bar{K}$, the algebraic closure of $K$ (see \cite[(A.8)]{W}).
The $\mu$-isomorphism induces a (right $H$-comodule) algebra isomorphism $(A/M)\otimes_{k}H\rightarrow (A/M)\otimes_{K}A$.
By extending this, we have a (right $H$-comodule) algebra isomorphism
\[ \bar{K}\otimes_{k}H\xrightarrow{\sim}\bar{K}\otimes_{K}A. \]
Therefore we have that (a) $\Leftrightarrow$ $\bar{K}\otimes_{K}A$ is reduced $\Leftrightarrow$ $A$ is a (classically) separable $K$-algebra.
(b) $\Leftrightarrow$ (c) is clear.
\end{proof}

If $k$ is an algebraically closed field in characteristic $0$, then
the above conditions are always satisfied.
But this fails in positive characteristic as was seen in Section \ref{sec-example}.

Then we generalize the Galois correspondence in Okugawa's sense.
Suppose that $k$ is algebraically closed and $H$ is reduced.
Let $\mathcal{F}$ be the set of intermediate AS $D$-module algebras $F$ of $L/K$ such that $AF\otimes_{F}AF$ is reduced.
When $L$ is a field, $\mathcal{F}$ is the set of intermediate $D$-module fields of $L/K$ over which $L$ is (classically) separable.
Let $\mathcal{G}$ be the set of all closed subgroups of $\bG(L/K)(k)$.
By Corollary \ref{cor} and Proposition \ref{prop-separable}, we have:
\begin{theorem}
\label{thm-Okugawa}
In the above situation,
$\mathcal{F}$ and $\mathcal{G}$ corresponds bijectively by
$\mathcal{F}\rightarrow\mathcal{G}$, $F\mapsto\bG(L/F)(k)$, and
$\mathcal{G}\rightarrow\mathcal{F}$, $\mathcal{H}\mapsto L^{\mathcal{H}}$.
\end{theorem}

\section*{Appendix}

Probably the following part is superfluous.
But I remain this because it includes a correction of an old version of my doctoral thesis.

Consider the Picard-Vessiot extension in Section \ref{sec-example}.
Then $L^{\bG(L/F)(k)}\neq F$ there.
But $F$ is the quotient field of $H/\langle l^{p}\rangle$-coinvariants of $A=K[f]$:
\[ A^{\mathrm{co}\,H/\langle l^{p}\rangle}=\{a\in A\;|\;\rho(a)-a\otimes 1\in A\otimes_{k}\langle l^{p}\rangle\}=K[f^{p}] \]
where $\rho : A\rightarrow A\otimes_{k}H$ is the right $H$-comodule structure of $A$.

Let $k$ be an arbitrary field and $D$ a cocommutative pointed $k$-Hopf algebra such that $D^{1}$ is a Birkhoff-Witt bialgebra.
Let $(L/K,A,H)$ be a Picard-Vessiot extension of AS $D$-module algebras (which is not necessarily finitely generated) such that $L^{D}=K^{D}=k$ and $F$ an intermediate AS $D$-module algebra corresponding to a Hopf ideal $I\subset H$.
Consider the following condition:
\begin{equation}
\label{FC}
\mbox{$H$ is a faithfully coflat $H/I$-comodule.}
\end{equation}
This means that $H\Box_{H/I}-$ is a faithfully exact functor where $\Box_{H/I}$ denotes the co-tensor product (see \cite{T2}).
Write $H':=H^{\mathrm{co}\,H/I}$.
If (\ref{FC}) is satisfied, then $H$ is a faithfully flat $H'$-module and $I=H(H')^{+}$ by \cite[Theorem 3]{T2}.
Especially this is the case when $I$ is a normal Hopf ideal.
\begin{theorem}
\label{thm-FC}
If the condition (\ref{FC}) is satisfied, then $F$ is the total quotient ring of $A^{\mathrm{co}\,H/I}$ ($=A\cap F$).
\end{theorem}
\begin{proof}
(This is a correction of the confused ``proof" of Proposition 3.5.7 (iii)
in the version 2 of \cite{Thesis}.)
We have that $A^{\mathrm{co}\,H/I}=A\cap F$ and the $\mu$-isomorphism $A\otimes_{k}H\xrightarrow{\sim}A\otimes_{K}A$ induces an isomorphism $A\otimes_{k}H'\xrightarrow{\sim}A\otimes_{K}(A\cap F)$ by \cite[Proposition 3.5.7 (ii)]{Thesis}.

Let $F'$ be the total quotient ring of $A\cap F$ realized in $L$.
Then $F'$ is an intermediate AS $D$-module algebra of $L/K$ which is included in $F$ by \cite[Corollary 3.3.8]{Thesis}.
Let $I'\subset H$ be the Hopf ideal corresponding to $F'$.
Then $I'\subset I$ since $F'\subset F$.
On the other hand, we have $H^{\mathrm{co}\,H/I'}\supset H'$ since $A\cap F'\supset A\cap F$ and $A\otimes_{k}H^{\mathrm{co}\,H/I'}\xrightarrow{\sim}A\otimes_{K}(A\cap F')$.
By the condition (\ref{FC}), $I=H(H')^{+}$.
Recall that the counit of $H$ is the restriction of $\mathrm{mult} : A\otimes_{K}A \rightarrow A$, $a\otimes b\mapsto ab$.
Thus,
\[ \begin{array}{l}
I=H(H')^{+}\subset H(H^{\mathrm{co}\,H/I'})^{+}\subset H\cap(A\cdot\kernel(A\otimes_{K}(A\cap F')\xrightarrow{\mathrm{mult}}A)\cdot A) \vspace{5pt} \\
\subset H\cap\kernel(L\otimes_{K}L\twoheadrightarrow L\otimes_{F'}L)=I'.
\end{array} \]
Therefore $I=I'$ and hence $F=F'$.
\end{proof}


\bibliographystyle{amsalpha}

\end{document}